\documentstyle[preprint,aps]{revtex}

\tightenlines

\begin{document}

\draft

\title{Integrable Systems and Harmonic Maps into Lie Groups}

\author{Fergus O'Dea}

\address{Department of Mathematics, National University of Ireland,
Galway, Ireland}

\date{\today}

\maketitle

\begin{abstract}

We study harmonic maps into Lie groups as a generalisation of the study of other well-known integrable systems, particularly the Toda and self-dual Chern Simons theories.

\end{abstract}

\section{Introduction}

The concept of integrability in infinite dimensions is not clear cut.  
In classical mechanics, a Hamiltonian system with
a $2n$-dimensional phase space is said to be
integrable  if it has
$n$ constants of the motion in involution (with vanishing Poisson brackets).
In systems of partial diferential equations, there
are infinitely many degrees of freedom, and there
is no straight-forward corresponding definition. 
Since the equations we are working with are static equations, not flow equations, integrals are not even defined.

There are other features of the systems that are
associated with being
`integrable'.  First, the equations are to some
degree soluble, meaning that explicit solutions
can be found and there exist general methods for
constructing solutions, which may be
superimposable in some extended sense. Alternatively, it may be possible
to find a large number of constants of the motion,
or the system may have the Painleve
 property.

The prototype of the integrable system is the Korteweg-De Vries (KdV) equation.  A solitary
wave is modelled by a soliton solution to the KdV equation:

\begin{equation}
$$4u_t-u_{xxx}-6uu_x=0$$\label{KdV}
\end{equation}
where $\kappa$ and $x_0$ are constants. The integrability of the KdV equation can be traced to
the  existence of a Lax pair, or obtaining the equation is the condition that the two differential operators
commute.  For KdV, the equation is the condition that
the two differential operators

\begin{equation}
$$L=\partial _x^2+u\ \ \ \ \ \ and\ \ \ \ \ \ M=\partial _t-\partial _x^3-{3 \over 2}u\partial _x-{3 \over 4}u_x$$\label{KdV}
\end{equation}

commute.

In this paper we look at the nature of relationships between some of
the better-known integrable systems. 
We focus on one of the most widely studied integrable
system,  the Toda model, and find a direct connection to the theory of harmonic maps into Lie
groups.

\section{Toda Theory}

The Toda field is a multicomponent field in two dimensions satisfying

\begin{equation}
$$\partial ^2\phi _i+\sum {k_{ij}}e^{\phi
_j}=0$$\label{toda}
\end{equation}
where $k_{ij}$ is the Cartan matrix of a
semisimple Lie algebra ${\cal g}$.  The Liouville equation is
obtained from  (\ref{toda}) by taking $k$ as the
Cartan matrix of the Lie algebra sl(2,{\bf C}).
The Toda field equations have a formulation in terms of a Lax pair.
First, we need some results from the theory of semisimple Lie algebra to describe this.   
If a Lie algebra $\cal L$
has a basis ${L^i}$ and structure constants $f^ij_k$ then

\begin{equation}
$$\left[ {L^i,L^j} \right]=\sum\limits_k {f_k^{ij}}L_k$$\label{struc}
\end{equation}
Associated with $\cal L$ is its set of roots $\Phi$.  These roots are $r$-dimensional vectors, 
where the rank $r$ of $\cal L$ is the maximal number of linearly independent commuting generators.
A set of simple roots $\Delta ={\alpha _i;\ i=1...r}$ is a subset of $\Phi$ such that
the difference of any two of its elements is not in $\Phi$.  Any root may be expressed

\begin{equation}
$$\alpha =\sum\limits_i {n_i\alpha _i}$$\label{}
\end{equation}
with the $n_i$ either all non-negative or all non-positive integers.  Hence $\Phi$ is divided
into two sets $\Phi^+$ and $\Phi^-$ containing positive and negative roots respectively.

The Cartan subalgebra is the maximal abelian subalgebra of $\cal L$ is called the Cartan
subalgebra. For $SU(n)$, these are just
the  diagonal matrices. A Chevalley basis of the Lie algebra is a choice of 
an element $H_i$ of the Cartan
for each simple root. The remaining
generators
$E_\alpha$ are labelled by roots.  Those labelled
by simple roots are for convenience denoted as
$E_i^\pm :=E_{\pm \alpha _i}$, which for our
purposes are just the generators for the
off-diagonal parts of the algebra.

In this basis the algebra takes the form

\begin{equation}
$$\matrix{\left[ {H_i,H_j} \right]=0\hfill\cr
  \left[ {H_i,E_j^\pm } \right]=\pm k_{ji}E_j^\pm \hfill\cr
  \left[ {E_i^+,E_j^-} \right]=\delta _{ij}H_j\hfill\cr}$$
 \label{chevalg}
\end{equation}

The upper and
lower triangular matrices are given by

\begin{equation}
$${\cal L}_\pm =\left\{ {H_i,E_\alpha :\ i=1...r,\
\alpha \in \Phi ^\pm } \right\}$$\label{}
\end{equation}

We are now ready to write down the Lax pair and the associated linear problem which
codify the Toda theory.

We follow Leznov and Saveliev \cite{LS} and
consider a connection whose components have values in different
subalgebras

\begin{equation}
$$A_{z,{\bar z}} \in {\cal L}_{+,-}
\end{equation}

If the connection is trivialisable, that is, it satisfies the linear equations

\begin{equation}
$$\partial _zg=A_zg,\ \ \ \ \ \ \ \partial _{\bar z}g=A_{\bar z}g$$
\end{equation}

then it automatically satisfies the flatness or integrability condition

\begin{equation}
$$\left[ {\partial _z+A_z,\partial _{\bar z}+A_{\bar z}} \right]=0$$
\label{flatness}
\end{equation}

An appropriate choice for $A$ is:

\begin{equation}
$$\matrix{A_z=\sum\limits_i {\left( {\partial _z\psi _-H_i+\alpha E_i^+} \right)}\hfill\cr
  A_{\bar z}=\sum\limits_i {\alpha e^{\beta \phi _i}}E_i^-\hfill\cr}$$\label{}
\end{equation}
where $\alpha, \beta$ are constants and $\psi,\phi$ are two $r$-component classical fields.
The algebra is straightforward, and in
characteristic fashion we get an
equation for each element of the
Chevalley basis that appears. Write out
(\ref{flatness}) explicitly, to get

\begin{equation}
$$\matrix{\ \ \ \partial _zA_{\bar z}-\partial
_{\bar z}A_z+\sum\limits_{ij} {\left( {\partial _+\psi
_i\alpha e^{\beta \phi _j}\left[ {H_i,E_j^-}
\right]+\alpha ^2e^{\beta \phi _j}\left[
{E_i^+,E_j^-} \right]} \right)}\hfill\cr
  =\sum\limits_i {\left( {\alpha \beta \partial _+\phi _ie^{\beta \phi _i}E_i^--\partial
_z\partial _{\bar z}\psi _iH_i} \right)}-\sum\limits_{ij} {\alpha \partial _+\psi _ie^{\beta
\phi _j}k_{ji}E_j^-}+\sum\limits_i {\alpha ^2}e^{\beta \phi _i}H_i\hfill\cr
  =0\hfill\cr}$$
\label{}
\end{equation}

which is satisfied by

\begin{equation}
$$\matrix{\psi _i=\beta \sum\limits_j {k_{ij}^{-1}}\phi _j\hfill\cr
  \partial _z\partial _{\bar z}\psi _i=\alpha ^2e^{\beta \phi _i}\hfill\cr}$$
\end{equation}

or written in the form of a single equation, the
Toda equation itself

\begin{equation}
$$\partial _z\partial _{\bar z}\phi _i={{\alpha ^2} \over \beta }\sum\limits_j
{k_{ij}}e^{\beta
\phi _j}$$\label{gentoda}
\end{equation}

\section{Conformal Affine Toda Theory}

Because of the appearance of the spectral parameter $\lambda$ in our Lax pair for the harmonic
map, we will find it is useful to study Lie algebras containing an affine parameter,
i.e., affine Lie algebras. The corresponding Toda
theories are the Affine Toda theories.  We will be
able to describe a relationship between 
harmonic maps into Lie groups and all of these models, most
interestingly with the so-called Conformal Affine
Toda theory.  This conformally invariant field
theory is based on the affine Lie algebra
$\hat{sl_2}$ which reduces under certain
circumstances to the Liouville theory and the
non-conformal sinh-Gordon theory. 

We first
describe the conformally invariant CAT theory \cite{BB}.
These models are
obtained from the usual Toda field theory by
adding two fields which transform in the correct
way under conformal transformations.   The
addition of these fields is facilitated by
constructing the affine Lie algebra
$\hat{sl_2}$.   This is the Lie algebra of
traceless 2x2 matrices with entries which are
Laurent polynomials in $\lambda$ (the loop
algebra  $\tilde{sl_2}$).  This algebra is
centrally extended as follows:

\begin{equation}
\mathop {sl'_2}\limits^{\wedge} =\mathop {sl_2}\limits^{\sim}\oplus
{\cal C}c \ .	\label{eq2}
\end{equation}
where 
\begin{equation}
\left[ {\hat X,\hat Y}
\right]_\wedge =\left[ {{\tilde X},{\tilde Y}} \right] _{\sim}+{1 \over {2i\pi
}}\oint {d\lambda \ tr\left[ {\partial _{\lambda} {\tilde X}\left( \lambda
\right)\cdot {\tilde Y}\left( \lambda  \right)} \right]} c\label{eq1}  
\end{equation} 
The affine Lie algebra $\hat{sl_2}$ is obtained by
adding the derivation $d=\lambda\frac{d}{d\lambda}$.  The algebra
can be decomposed as
\begin{equation}
\mathop {sl'_2}\limits^\wedge ={\cal N}_-\oplus {\cal{H}}\oplus 
{\cal N}_+\label{eq3}
\end{equation}
where ${\cal N}_-$,{\em}\ ${\cal N}+$ are lower and upper triangular matrices
respectively, and $\cal{H}$, the Cartan sub-algebra, is spanned by the elements
H, c and d.  We write in co-ordinates $\Phi: {\cal C} \to {\cal H}$

\begin{equation}
\Phi ={1 \over 2}\phi H+\eta d+{1 \over 2}\xi c\label{eq4}
\end{equation}
where $H, d, c$ generate ${\cal H}$.

We follow the construction of Toda field theory shown above, but using complex co-ordinates
and a slightly different connection. The Lax pair
$(\partial_z+{\cal{A}}_z, \partial_{\bar z}+{\cal{A}}_{\bar z})$ is written as
\begin{eqnarray}
A_z&=&\partial _{z}\Phi +e^\Phi (E_++\lambda E_-)e^{-\Phi } \\
A_{\bar
z}&=&-\partial _{\bar z}\Phi +e^{-\Phi }(E_-+\lambda
^{-1}E_+)e^\Phi \ . \label{eq5} 
\label{laxp}
\end{eqnarray}

The zero
curvature condition (\ref{flatness}) now gives us,
after a little effort, the following set of
equations
\begin{eqnarray}
\partial _z\partial _{\bar z}\phi \label{eqp1}
&=&e^{2\phi }-e^{2\eta -2\phi } \\ \partial _z\partial _{\bar z}\eta
&=&0 \\ \partial _z\partial _{\bar z}\xi &=&e^{2\eta -2\phi } \ . \label{eq6}
\end{eqnarray}

This system of three equations is the Conformal Affine Toda theory.  The reduction to the sinh-Gordon and Liouville 
theories are realised 
by the limits $\eta\rightarrow 0 $ and
$\eta\rightarrow-\infty$ respectively.

The set of equations
(\ref{eqp1})-(\ref{eq6}) are conformally
invariant, as is the
 reduction to the Liouville equation.

\section{Wess-Zumino-Witten models and Liouville theory}

O'Raifeartaigh et al \cite{FWBFOR} have shown that 
Liouville theory can be regarded as a reduced SL(2,R) Wess-Zumino-Witten
theory. Their reduction uses a decomposition into local fields.

  The energy or action for a group-valued field $g$ is

\begin{equation}  
$$S(g)=-{k \over {8\pi }}\int {d^2}z Tr\left[ {(g^{-1}\partial _z g)(g^{-1}\partial
_{\bar z}g)}
\right]+{k \over {12\pi }}\int\limits_{B^3} {Tr\left[ {(g^{-1}dg)^3} \right]}$$
\end{equation}    

Any connected semi-simple real Lie group $G$ admits a Gauss decomposition

\begin{equation}  
$$G=XYZ$$    \end{equation}    
where $Y$ is the direct product $Y=A\otimes K$ of a simply-connected abelian group $A$
and a connected semisimple compact group $K$, and the groups $X$ and $Z$ are simply
connected and nilpotent; two arbitrary decompositions are connected by an automorphism of $G$. \cite{BR}

In the case of $SL(2,R)$, there exists such a decomposition for regular g which valid in a
neighbourhood of the identity, as follows:

\begin{equation}
$$g=ABC$$
\label{gauss}
\end{equation}

where

\begin{eqnarray}
A=\left( {\matrix{1&x\cr
0&1\cr
}} \right)=\exp (xE_+) \\ 
  C=\left( {\matrix{1&0\cr
y&1\cr
}} \right)=\exp (yE_-) \\ 
  B=\left( {\matrix{{\exp ({\textstyle{1 \over 2}}\phi )}&0\cr
0&{\exp (-{\textstyle{1 \over 2}}\phi )}\cr
}} \right)=\exp ({\textstyle{1 \over 2}}\phi H)
\end{eqnarray}

The Wess-Zumino energy for the product of three matrices $A,B,C$
can be written as the sum of the energies for the actions for $A,B$ and $C$ plus another
term:

\begin{eqnarray}
S(ABC)=S(A)+S(B)+S(C) -{k \over {4\Pi }}\int {d^2\xi \ Tr}\ [ \left( {A^{-1}\partial _zA}
\right)\left( {\partial _{\bar z}B} \right)B^{-1} \\ +\left( {B^{-1}\partial _{\bar
z}B}
\right)\left( {\partial _zC} \right)C^{-1}+\left( {A^{-1}\partial _{\bar z}A}
\right)B\left( {\partial _zC} \right)C^{-1}B^{-1} ]
\end{eqnarray}

 Using our parametrization, the energy for the Wess-Zumino model takes the following
local form:

\begin{eqnarray}
\partial _z\partial _{\bar z}\phi +2e^{-\phi }\partial _zy\partial _{\bar z}x=0 \\
  \partial _z(\partial _{\bar z}xe^{-\phi })=\partial _{\bar z}(\partial _zye^{-\phi
})=0
 \end{eqnarray} 

The equations of motion for the Wess-Zumino model can be derived 
in these co-ordinates as

\begin{eqnarray}
\partial _z\partial _{\bar z}\phi +2e^{-\phi }\partial _zy\partial _{\bar z}x=0 \\
  \partial _z(\partial _{\bar z}xe^{-\phi })=\partial _{\bar z}(\partial _zye^{-\phi
})=0
\end{eqnarray}

Consider special solutions

\begin{equation}  
$$\partial _{\bar z}x=\nu e^\phi \ \ \ \ \ \ \ \partial_z y=\mu e^\phi $$
\label{orc}
 \end{equation}    
where $\mu$, $\nu$ are arbitrary constants.  Then the system
reduces to the Liouville system

\begin{equation}  
$$\partial_z \partial_{\bar z} \phi + 2\mu \nu e^\phi =0$$
\end{equation}

\section {Self-dual Chern Simons and Toda theory}

In \cite{D}, Gerald Dunne defined the self-dual
Chern-Simons equations over $SU(N)$ as:

\begin{eqnarray}
\partial _zA_{\bar z}-\partial _{\bar z}A_z+\left[ {A_{\bar z},A_z} \right]&=&{2 \over
\kappa }\left[ {\Psi ^{\dagger},\Psi } \right] \label{38} \\
  \partial _{\bar z}\Psi +\left[ {A_{\bar z},\Psi } \right]&=&0 \ . \label{39}
\end{eqnarray}

He showed that, with for a diagonal $A$ and upper triangular $\Psi$, 

\begin{eqnarray}
A_z=\sum\limits_{}^{} {A_i^\alpha }H_\alpha \\
  \Psi =\sum {\psi ^\alpha }E_\alpha 
\end{eqnarray}

these equations combine
to become those of the Toda model:

\begin{eqnarray}
\partial ^2\phi _\alpha =-{2 \over \kappa }K_{\alpha \beta }\phi _\beta \ .
\end{eqnarray}

where $\ln \phi _\alpha \equiv \left| {\psi ^\alpha } \right|^2$.   He also showed it is
possible to make a gauge transformation
$u^{-1}$ which combines the self-dual Chern-Simons equations into a
single equation:

\begin{equation}
$$\partial _{\bar z}\chi =\left[ {\chi ^{\dagger},\chi } \right]$$
\label{40}
\end{equation}
where 

\begin{equation}
$$\chi =\sqrt {{2 \over \kappa }}u\Psi u^{-1}$$
\label{44}
\end{equation}

Define

\begin{eqnarray}
\tilde A_z\equiv A_z-\sqrt {{2 \over \kappa }}\Psi \\
  \tilde A_{\bar z}\equiv A_{\bar z}+\sqrt {{2 \over \kappa }}\Psi ^{\dagger}
\end{eqnarray}

Then the self-dual Chern Simons equations imply
that $\tilde A$ is flat.  Trivialising  $\tilde
A$ as  $\tilde A=u^{-1}du$, and using the $\chi$
defined above,
we find that

\begin{eqnarray}
\partial _zA_{\bar z}-\partial _{\bar z}A_z+\left[ {A_{\bar z},A_z} \right]-{2 \over
\kappa }\left[ {\Psi ^{\dagger},\Psi } \right]=
  g^{-1}\left( {\partial _{\bar z}\chi +\partial _z\chi ^{\dagger}-2\left[ {\chi
^{\dagger},\chi } \right]} \right)g \ .
\end{eqnarray}
This shows that the equations (\ref{38}) and (\ref{39}) are
equivalent to the single equation (\ref{40}).

This equation may now be written as the harmonic
map equation 

\begin{eqnarray}
\partial _z\left( {h^{-1}\partial _{\bar z}h} \right)+\partial _{\bar z}\left(
{h^{-1}\partial _zh} \right)=0 \ ,
\end{eqnarray}
where $h \in SU(N)$ is related to $\chi$ as

\begin{eqnarray}
h^{-1}\partial _zh=2\chi \ .
\end{eqnarray}
All SU(2) finite action harmonic
maps have the form 

\begin{eqnarray}
h=-h_0\left( {2p-1} \right)
\end{eqnarray}
where $p$ is a holomorphic projection valued map

\begin{eqnarray}
\left( {1-p} \right)\partial _zp=0 \ .
\end{eqnarray}
With this condition, we can write a general $p$ in the defining representation for SU(2) as

\begin{eqnarray}  
p={{MM^\dagger} \over {M^\dagger M}} 
\end{eqnarray}     where 

\begin{eqnarray}  
M=\left( {\matrix{1\cr
{f(z)}\cr
}} \right)
   \end{eqnarray}
for arbitrary $f(z)$.  It can easily be checked that $p$ satisfies the correct projectivity, hermiticity and holomorphicity
conditions.

Hence, we find that

\begin{equation}  
$$p={1 \over {1+f\bar f}}\left( {\matrix{1&{\bar
f}\cr f&{f\bar f}\cr
}} \right)$$    \end{equation}    
 
Explicitly, the $\chi$ become

\begin{eqnarray} 
{\textstyle{1 \over 2}}h^{-1}\partial
_zh={\textstyle{1 \over 2}}(2p-1)\cdot 2\partial
_zp=\partial _zp \\
  ={{f\partial _z\bar f} \over {\left( {1+f\bar f}
\right)^2}}\left( {\matrix{{-1}&{{1 \over f}}\cr
{-f}&1\cr }} \right) \ .
    \end{eqnarray}    
The corresponding commutator is

\begin{equation}
$$\left[ {\chi ,\chi ^{\dagger}} \right]={{\partial _z\bar f\partial _{\bar z}f} \over
{\left( {1+f\bar f} \right)^3}}\left( {\matrix{{1-f\bar f}&{2\bar f}\cr {2f}&{-1+f\bar f}\cr
}} \right)$$
\end{equation}

Dunne showed that this can be diagonalised by $u\in SU(2)$

\begin{equation}  
$$u={1 \over {\sqrt {1+f\bar f}}}\left(
{\matrix{1&{-\bar f}\cr f&1\cr
}} \right) \ , $$   
\label{u} \end{equation} 

so that

\begin{equation}  
$$u^{-1}\left[ {\chi ,\chi ^{\dagger}} \right]u={{-\partial _z\bar f\partial
_{\bar z}f} \over {\left( {1+f\bar f} \right)^2}}\left(
{\matrix{1&0\cr 0&{-1}\cr }} \right) . $$   
\end{equation}

Using (\ref{44}), we see that $\ln \phi _\alpha \equiv \left| {\psi ^\alpha } \right|^2$ is the quantity
$${{-\partial _z\bar f\partial _{\bar z}f} \over {\left( {1+f\bar f} \right)^2}}$$
where $\phi$ is the Toda field.  It above can also be written as 

\begin{equation}  
$$\partial _z\partial _{\bar z} \ln (1+f\bar f)=\partial _z\partial _{\bar z} \ln \
\det(M^\dagger M)$$   
\label{md}
 \end{equation}    

This is the form of the general solution of the
classical Liouville equation  \cite{L}.

\section{Harmonic Maps and Toda Systems}

We seek a broader understanding of the relationship of harmonic maps 
to other integrable systems.
 In the process we construct a general
scheme for reducing certain Harmonic maps into Lie groups
 to the Toda system and examine the relationships with other 
 models discussed in the previous chapter.

We look mainly at harmonic maps
into SL(2,R) or equivalently, SU(1,1), a choice which allows us to compare
 our constructions to both the work of Dunne \cite{D} and
O'Raifertaigh et. al.  \cite{FWBFOR}. Using  analysis described
by Uhlenbeck in
\cite{U} for harmonic maps into unitary groups, we
find that the work of Dunne can be used to find
solutions for the Liouville system. The second
section of this chapter employs the
field-theoretic approach of O'Raifertaigh et. al.
 using a decomposition of the group
SU(1,1) to demonstrate an alternative reduction to the
Toda model.  

We then construct a new framework which puts these reductions 
in a single context.  A requirement that the fields in the Lax pair satisfy certain
requirements in their Lie algebraic structure is found to be behind the conditions and
{\it Ansatze} used in \cite{FWBFOR}, \cite{D} and \cite{G}, and others.

There exists a relatively straightforward way to associate a harmonic map and
the corresponding harmonic map (chiral model) equations with a CAT theory. 
This method also holds for reducing the non-affine algebra (and the corresponding
system of equations) to the simple Toda theory.

The correspondence relies on the following
proposition, first proved in \cite{U2}:

{\bf Proposition 1} {\it If $A=A_z
dz+A_{\bar{z}}d \bar {z}$ and
$B=B_z dz+B_{\bar{z}}d\bar {z}$ satisfy}

\begin{equation}
\left[ {\partial _{\bar z}+A_{\bar z}+\lambda B_{\bar z},\partial
_z+A_z+\lambda ^{-1}B_z} \right]=0\ ,\label{eq8}
\end{equation}
{\it then there exists $u$ with}
\begin{equation}
\tilde A=u^{-1}(Au+du),\ \tilde B=u^{-1}Bu\label{eq8a}
\end{equation}

{\it and}
\begin{equation}
\left[ {\partial_{\bar z}+(1-\lambda )\tilde A_{\bar z},\partial _z+(1-\lambda
^{-1})\tilde A_z} \right]=0\ .\label{eq9}
\end{equation}

{\it Furthermore, $2\tilde{A}=s^{-1}ds$ where $s$ is harmonic.}

{\it Proof:}  The necessary transformations correspond to trivialising
(\ref{eq8}) at $\lambda=1$, with the gauge transformation $u$. 
From \cite{U}, since ${\bf C}$ is simply
connected, we find that $s$ is a harmonic map where
$\tilde{A}=s^{-1}ds/2$.

Notice that (\ref{laxp}) can be written in the correct format to allow
use of Proposition 1, where
\begin{eqnarray}
&A_z&=\partial _z\Phi +e^{\Phi}E_+e^{-\Phi}\\   
&A_z&=-\partial _{\bar z}\Phi+e^{\Phi}E_-e^{-\Phi} \\ 
&B_z&=e^{\Phi}E_-e^{-Phi} \\  
&B_{\bar z}&=e^{\Phi}E_+e^{-\Phi} \ .
\label{eq10}
\end{eqnarray}

We write the Lax pair of our harmonic map as

\begin{equation}
$$\left( {\partial _{\bar z}-(1-\lambda )\ u^{-1}e^\Phi E_+e^{-\Phi }u,\partial _z-(1-\lambda ^{-1})u^{-1}e^\Phi E_-e^{-\Phi }u} \right) \ ,$$
\end{equation}
an equation which holds for all order of $\lambda$.

\section{A condition for the equivalence of the harmonic map and
Conformal Affine Toda equations}

Not all harmonic maps are Toda maps.
We now show that Toda maps, and more generally, CAT maps, are a
subset of harmonic maps into the appropriate Lie
group.  

As a brief introduction on the way to a more general result for all
semisimple Lie groups, we now describe a condition under
which a harmonic map into SL(2,C)' (the subgroup of
constant loops in the centrally extended loop
group) gives rise to the CAT system of
equations.   

{\bf Theorem 2}\ {\it If $s$ is a harmonic map into
SL(2,C), then the harmonic map equations for $s$
are gauge equivalent to the Conformal Affine Toda
equations (\ref {eq6}) for
$\phi, \eta, \xi$ if the fields $s^{-1} \partial
_z s, s^{-1} \partial_{\bar z} s$ can be
simultaneously diagonalised into lower and
upper triangular matrices respectively.}

{\it Proof:}  Let $\tilde{A}_{z} =s^{-1}\partial_z
s/2, \tilde{A}_{\bar z} =s^{-1}\partial_{\bar z} s/2$. Thus we can write, for
some $\theta \in \hat {sl_2}'$,

\begin{equation}
\tilde A_z=e^{-\theta}fE_+e^{\theta},\ \tilde A_z=e^{-\theta}f'E_-e^{\theta}\label{eq13}
\end{equation}

By proposition 1, the harmonic map equations
for
$s$ are described by the vanishing of the Lax pair
in (\ref{eq9})  for all
$\lambda \in C^*$.  In terms of the connection
given above, this provides the following
identites as coefficients for $\lambda^{-1},
\lambda$ respectively: 
\begin{eqnarray}
-\partial _{\bar z}fE_+-f\left[
{\partial _{\bar z}\theta,E_+} \right]-2f'fH=0 \\ \partial _zf'E_++f'\left[
{\partial _z \theta,E_-} \right]-2f'fH=0
\label{eq15}
\end{eqnarray}

Writing $\partial_z\theta, \partial_{\bar z}\theta$ in terms of its algebraic decomposition

\begin{eqnarray}
\partial _z\theta=\alpha E_++\beta
H+\gamma E_-+\delta c \\ \partial _{\bar
z}\theta=\alpha 'E_++\beta 'H+\gamma 'E_-+\delta
'c 
\label{eq16}
\end{eqnarray}
we see that the following relationships are
specified:
\begin{eqnarray}
\partial _zf'+\beta f'=0\ \ \ \ \ \ \ \ \ \
\ \alpha f-f'f=0  \\ \partial _{\bar z}f-\beta 'f=0\ \ \ \ \ \ \ \ \ \ \ \gamma
'f'-f'f=0 \label{eq17}
\end{eqnarray}

Using this information, we write
\begin{eqnarray}
\partial _z \theta=f'E_+-\partial _z\left( {ln\ f'} \right)H+\gamma E_-+\delta c \\ \partial _{\bar z}\theta=fE_-+\partial _{\bar z}\left( {ln\ f} \right)H+\alpha
'E_++\delta 'c \ . \label{eq18}
\end{eqnarray}
Apply Proposition 1 in exactly the reverse
manner using the inverse of the gauge transformation given by the
trivialisation at the identity.  In other words, transform to

\begin{eqnarray}
A=u\tilde Au^{-1}-du\cdot u^{-1}  \\ B=u\tilde Bu=-u\tilde Au^{-1}  
\label{eq20}
\end{eqnarray}
where here we have $u=e^\theta$.  Define the
quantities
\begin{eqnarray}
\phi _1&=ln\ f;\ \ \ \
\ \eta _1&=\ ln\ \gamma +\phi _1  \\ \phi
_2&=ln\ f';\ \ \ \ \eta _2&=\ ln\ \alpha '
+\phi _2\ \ .  \label{eq19}
\end{eqnarray} We find the connection becomes
\begin{eqnarray}
&\left( {\partial _{\bar z}+A_{\bar
z}+\lambda B_{\bar z},\partial _z+A_z+\lambda ^{-1}B_z} \right) \\ &=\left( {\partial _{\bar z}-\partial _{\bar z}\phi _1H-e^{\eta _1-\phi
_1}E_+-\delta 'c-\lambda e^{\phi _1}E_-,\ \ \partial _z+\partial _z\phi
_2H-e^{\eta _2-\phi _2}E_--\delta c-\lambda ^{-1}e^{\phi _2}E_+} \right) \ .   \label{eq21}
\end{eqnarray} The curvature is given by \begin{eqnarray}
&\left( {-\ \partial
_z\partial _{\bar z}\left( {\phi _1+\phi _2} \right)+2e^{\phi _1+\phi
_2}-2e^{\eta _1+\eta _2-\phi _1-\phi _2}} \right)H \\ &+\left( {-\partial _ze^{\eta
_1-\phi _1}+\lambda ^{-1}\partial _{\bar z}\phi _1e^{\phi _2}-\lambda
^{-1}\partial _{\bar z}\phi _2e^{\phi _2}-\partial _z\phi _2e^{\eta _1-\phi
_1}} \right)E_+ \\ &+\left( {\partial _{\bar z}e^{\eta _2-\phi _2}+\lambda
^{}\partial _z\phi _2e^{\phi _1}-\lambda ^{-1}\partial _z\phi _1e^{\phi
_1}+\partial _{\bar z}\phi _1e^{\eta _2-\phi _2}} \right)E_- \\ &+\left(
{-\partial _z\delta '+\partial _{\bar z}\delta +e^{\phi _1+\phi _2}} \right)c \ .
\label{eq22}
\end{eqnarray}

We now introduce the variables $\xi_1,\xi_2$ where
$\delta ' =\partial_{\bar{z}}\xi_1, \delta_2=\partial_z\xi_2$.  We find that,
by writing
\begin{eqnarray}
\eta =  \eta _1+\eta _2;\ \ \ \phi =\phi_1+\phi _2 ; \\ \xi =  \xi _1+\xi _2+\phi
\label{eq23}
\end{eqnarray} 

this system is equivalent to the CAT set:

\begin{eqnarray}
\partial _z\partial _{\bar z}\phi \label{eqp1}
&=&e^{2\phi }-e^{2\eta -2\phi } \\ \partial _z\partial _{\bar z}\eta \label{eqp2}
&=&0 \\ \partial _z\partial _{\bar z}\xi &=&e^{2\eta -2\phi } \ . \label{eq6}
\end{eqnarray}

Note that a reduction of our result to the
unextended sl(2) algebra will lead to a subset of
the CAT equations, namely equations (\ref{eqp1}) and (\ref{eqp2}). 
Although equation (\ref{eq6}) seems somewhat superfluous
to requirements, in that it can be safely omitted
while preserving the conformal invariance, the
central extension plays an important role by
including the spectral parameter into the
algebra.  Also, the conformal invariance
of the CAT theory can be described by the invariance of
the harmonic map equations for the system 
described above under the transformation 
$f \to {\tilde f}=g(z)f$ .

\section{Algebraic Reduction of Harmonic Maps to Toda Systems}

In the last section 
we reduced the harmonic map
equations  to the Conformal Affine Toda systems for SL(2).   In this section, we will
investigate harmonic maps into SU(1,1) and see
how we can identify
these with solutions of the Toda equations.  Later
we will generalise this to include all semi-simple
Lie groups.

Let us first construct the analog of Uhlenbeck's
uniton for the case of a harmonic map into
SU(1,1). We will begin by generalising the uniton
construction of \cite{U}  to the indefinite
$U(q,N-q)$ form.  Again, we consider maps into the Grassmannian of $k$ planes
which is a geodesic submanifold.

{\bf Definition}:  An $n,q$-uniton is a harmonic map $s:\Omega \rightarrow U(q,N-q)$ which has an
extended solution

\begin{equation}  
$$E_\lambda :{\bf C}^*\times \Omega \to G=GL(N,C)$$    \end{equation}    
with

\begin{eqnarray}
&(a)&\ E_\lambda =\sum\nolimits_{\alpha
=0}^n {T_\alpha }\lambda ^{\alpha \ }\ for\
T_\alpha :\Omega \to gl(N,C) \\
  &(b)&\ E_1=I \\
  &(c)&\ E_{-1}=Qs^{-1}\ for\ Q\in SU(q,N-q)\
constant \\
  &(d)&\ (E_{\bar \lambda })^*J=J(E_{\lambda
^{-1}})^{-1}
\end{eqnarray}

and $J$ is the automorphism which defines the real form for $U(q,N-q)$

{\bf Proposition 3}: {\it $s:\Omega \rightarrow
SU(q,N-q)$ is a 1-q-uniton (a holomorphic map) if 
$s=Q(2p-1)$ for constant $Q\in SU(q,N-q)$, where
$p$ satisfies $p^*J=Jp$, $p^2=p$, and
$(1-p)\partial _{\bar z}p=0$ }

{\bf Proof}: With $E_\lambda =T_0+\lambda T_1$
 and $E_1=T_0+T_1=I$, we have
$E_\lambda =T_0+\lambda (I-T_0)$.

The reality condition (d) above tells us that

\begin{eqnarray}
(I-T_0^*)JT_0&=&0 \\
  T_0^*JT_0+(I-T_0^*)J(I-T_0)&=&J
\end{eqnarray}

Combining these equations gives $T_0^*J=JT_0$ and $T_0^2=T_0$.
We can identify $T_0$ with $p$.  As in the previous chapter
the necessary condition for $E_\lambda$ to 
be an extended solution for a harmonic map as
is that $(1-p)\partial _{\bar z}p=0$

With this condition, we can write a general $p$ in the defining representation for SU(1,1) as

\begin{equation}  
$$p=M(M^*JM)^{-1}M^*J$$    \end{equation}    
where, as before, 

\begin{equation}  
$$M=\left( {\matrix{1\cr
{f(x_-)}\cr
}} \right)$$    \end{equation}    
for arbitrary $f(x_-)$. 

Hence, we find that, in this case

\begin{equation}  
$$p={1 \over {1-f\bar f}}\left( {\matrix{1&{-\bar f}\cr
f&{-f\bar f}\cr
}} \right)$$    \end{equation}    
 
The simplest 
solutions of the harmonic map equation are now given by 
the holomorphic maps:

\begin{eqnarray} 
{\textstyle{1 \over 2}}g^{-1}\partial _zg={\textstyle{1 \over 2}}(2p-1)\cdot 2\partial
_zp=\partial _zp \\
  ={{-f\partial _z\bar f} \over {\left( {1-f\bar f} \right)^2}}\left( {\matrix{{-1}&{{1
\over f}}\cr {-f}&1\cr
}} \right)
    \end{eqnarray}

We can now find a $u\in SU(1,1)$

\begin{equation}  
$$u={1 \over {\sqrt {1-f\bar f}}}\left( {\matrix{1&{\bar f}\cr
f&1\cr
}} \right)$$   
\label{u} \end{equation} 

such that

\begin{equation}  
$$u^{-1}[A_+,A_-]u={{-\partial _z\bar f\partial _{\bar z}f} \over {\left( {1-f\bar f} \right)^2}}\left( {\matrix{1&0\cr
0&{-1}\cr
}} \right)$$    \end{equation}

This diagonalisation is a critical step, and    
 as was shown in the previous chapter for SU(2), the magnitude of this quantity should
be a solution of the appropriate Toda equation. It can be rewritten 

\begin{equation}  
$$u={1 \over {\sqrt {1-f\bar f}}}\left( {\matrix{1&{\bar f}\cr
f&1\cr
}} \right)$$   
\label{u} \end{equation}    
   
and we find that the magnitude of this expression for the diagonalised commutator above is just

\begin{equation}  
$$\ln \phi = \partial _z\partial _{\bar z}\ln (1-f\bar f)=\partial _z\partial _{\bar z} \ln
\
\det(M^*JM) \ .$$   
\label{md}
 \end{equation}

\section{Harmonic Maps and Liouville Theory}

Let us now examine the reduction of the SU(1,1) harmonic map
to the Toda system using a decomposition into local fields.

Recall that the energy or action for $g: \Omega \to SU(1,1)$ is

\begin{equation}  
$$S(g)=-{k \over {8\pi }}\int {d^2}z Tr\left[ {(g^{-1}\partial _z
g)(g^{-1}\partial _{\bar z}g)}
\right]$$
\end{equation}

As in the previous chapter, where we saw 
a similar reduction to
analyse Wess-Zumino-Witten theory,  we use a Gauss decomposition of
$SU(1,1)$.  For regular $g$ in a neighbourhood of the identity, we write:

\begin{equation}
$$g=ABC$$
\label{gauss}
\end{equation}

where

\begin{eqnarray}
A=\left( {\matrix{1&{ix}\cr
0&1\cr
}} \right)=\exp (ixE_+) \\
  C=\left( {\matrix{1&0\cr
{iy}&1\cr
}} \right)=\exp (iyE_-) \\
  B=\left( {\matrix{{\exp ({\textstyle{1 \over 2}}\phi )}&0\cr
0&{\exp (-{\textstyle{1 \over 2}}\phi )}\cr
}} \right)=\exp ({\textstyle{1 \over 2}}\phi H)    
\end{eqnarray}

With this parametrization, the energy for the harmonic map can be written
in terms of three real fields as

\begin{equation}  
$$S(g)=S(x,y,\phi )=-{k \over {8\pi }}\int {d^2}z \left[ {{\textstyle{1 \over 2}}\partial
_z\phi \partial _{\bar z}\phi -e^{-\phi }(\partial _zx\partial _{\bar z}y+\partial
_zy\partial _{\bar z}x)} \right]$$    \end{equation}     A local form of the equations of
motion for the harmonic map can now be derived from this:

\begin{eqnarray}
\partial _z\partial _{\bar z}\phi -e^{-\phi }(\partial _zx\partial _{\bar z}y+\partial
_zy\partial _{\bar z}x)=0 \\
  \partial _z(\partial _{\bar z}xe^{-\phi })+\partial _{\bar z}(\partial _zxe^{-\phi })=0
\\
  \partial _z(\partial _{\bar z}ye^{-\phi })+\partial _{\bar z}(\partial _zye^{-\phi
})=0 \ .
\end{eqnarray}

We consider the solution for the above:

\begin{equation}  
$$\partial _zx=f(z)e^\phi \ , \ \ \ \ \ \partial _zy=g(z)e^\phi \ .$$ 
\label{orc}
 \end{equation} 
   
Then the system 
reduces to the following form of the Liouville system:

\begin{eqnarray}
\partial _z\partial _{\bar z}\phi +M\left( z,{\bar z} \right)e^\phi =0 \\
  f\left( z \right)\bar g\left( {\bar z} \right)-\bar f\left( {\bar z} \right)g\left( z
\right)-M=0 \ .
\end{eqnarray}

\section{A useful result}

We have seen the reduction of the harmonic maps into Lie groups
to the Toda model carried out in a number of different ways, and there are others, such
as \cite{G}, which we have not alluded to. We now ask ourselves if there is any underlying
link in the manner in which these reductions are carried out.

We find that these reductions can be observed to 
be different strategies for taking advantage of the following general result:

{\bf Result 4} {\it Any pencil of connections $(A_z+\lambda A_z ', A_{\bar z}+\lambda
^{-1}A_{\bar z} ')$ gives rise to the Toda system of equations for the corresponding
semi-simple Lie group if the $A_z, A_{\bar z}$ are upper or lower triangular respectively
 and the $A_z ', A_{\bar z} '$
are off-diagonal.}

{\it Proof:}

We will begin with the most general case and reduce using our stated condition as
needed.  We see that most of the work is accomplished by the algebraic structure and the
parameter $\lambda$. 

We begin with the stated condition that

\begin{equation}  
$$\left[ {\partial _z+A_z+\lambda A'_z,\partial _{\bar z}+A_{\bar z}+\lambda ^{-1}A'_{\bar z}} \right]=0$$
\label{beg}
\end{equation}    

where $A$ and $A '$ take
values in the Chevalley basis of the Lie algebra, i.e.

\begin{equation}
$$\matrix{A_z=\sum\limits_{\alpha \in \Phi ^+} {g_+^\alpha }E_\alpha +\sum\limits_{\alpha
\in \Phi ^+} {f_+^\alpha }E_{-\alpha }+\sum\limits_{\alpha \in \Phi ^+} {h_+^\alpha
}H_\alpha \hfill\cr
  A_{\bar z}=\sum\limits_{\alpha \in \Phi ^+} {g_-^\alpha }E_\alpha +\sum\limits_{\alpha
\in
\Phi ^+} {f_-^\alpha }E_{-\alpha }+\sum\limits_{\alpha \in \Phi ^+} {h_-^\alpha }H_\alpha
\hfill\cr}$$
\label{63}
\end{equation}

and

\begin{equation}
$$\matrix{A'_z=\sum\limits_{\alpha \in \Phi ^+} {g'^\alpha_+ }E_\alpha
+\sum\limits_{\alpha \in \Phi ^+} {f'^\alpha_+ }E_{-\alpha }+\sum\limits_{\alpha
\in \Phi ^+} {h'^\alpha_+ }H_\alpha \hfill \cr
  A'_{\bar z}=\sum\limits_{\alpha \in \Phi ^+} {g'^\alpha_- E_\alpha }+\sum\limits_{\alpha
\in \Phi ^+} {f'^\alpha_- }E_{-\alpha }+\sum\limits_{\alpha \in \Phi ^+}
{h'^\alpha_- }H_\alpha
\hfill \cr}$$       
\end{equation}
\ 

Here  $\Phi ^+$ denotes the set of positive roots, and $H_\alpha, E_{\pm \alpha} $
are the Cartan subalgebra and step generators of the Chevalley basis respectively.

In terms of the first power of of $\lambda$ and the basis of the Lie algebra,
 the set of
equations reduces to:

\begin{eqnarray}
-\partial _{\bar z}g'^\alpha_+ -\sum\limits_\beta  {K_{\beta \alpha }}h_-^\beta g'^\alpha_+ +\sum\limits_\beta  {K_{\beta \alpha }}h'^\beta_+ g_-^\alpha =0 \\
  -\partial _{\bar z}f'^\alpha_+ +\sum\limits_\beta  {K_{\beta \alpha }}h_-^\beta f'^\alpha_+ -\sum\limits_\beta  {K_{\beta \alpha }}h'^\beta_+ f_-^\alpha =0 \\
  \partial _-h'^\alpha_+ +g'^\alpha_+ f_-^\alpha -f'^\alpha_+ g_-^\alpha =0
\label{67}
\end{eqnarray}
and their dual, where $K$ is the classical Cartan matrix for the Lie algebra.
In the case where $A'$ is off-diagonal, all the $h'$ vanish, so we find:

\begin{eqnarray}
\partial _z h_ - ^\alpha   - \partial _{\bar z} h_ + ^\alpha   + f_ - ^\alpha  g_ + ^\alpha   - f_ + ^\alpha  g_ - ^\alpha   + f'^\alpha_-   g'^\alpha_+    - f'^\alpha_+   g'^\alpha_-   = 0
\end{eqnarray}
The part of (\ref{beg}) independent of $\lambda$ contributes the following (and its dual:

 \begin{equation}
\partial _z g_ - ^\alpha   - \partial _{\bar z} g_ + ^\alpha   - \sum\limits_\beta  {K_{\beta \alpha } } h_ - ^\beta  g_ + ^\alpha   + \sum\limits_\beta  {K_{\beta \alpha } } h_ + ^\beta  g_ - ^\alpha   - \sum\limits_\beta  {K_{\beta \alpha } } h'^\beta_ -   g'^\alpha_ +    + \sum\limits_\beta  {K_{\beta \alpha } } h'^\beta_ +   g'^\alpha_ -    = 0
\end{equation}
as well as

 \begin{equation}
\partial _z h_ - ^\alpha   - \partial _{\bar z} h_ + ^\alpha   + f_ - ^\alpha  g_ + ^\alpha   - f_ + ^\alpha  g_ - ^\alpha   + f'^\alpha_ -   g'^\alpha_ +    - f'^\alpha_ +   g'^\alpha_ -   = 0
 \end{equation}

When $A$ is triangular, $g_+ ^\alpha =0$  and when $A'$ is
off-diagonal, $h'=0$ and we find $g'^\alpha_-=0$.  
Equation (\ref{67}) and its dual now give
us:

 \begin{eqnarray}
\partial _z\left( {\ln g_-^\alpha } \right)=-\sum\limits_\beta  {K_{\beta \alpha
}}h'^\beta_+ \\
  \partial _z\left( {\ln f_+^\alpha } \right)=\sum\limits_\beta  {K_{\beta \alpha
}}h'^\beta_-
\end{eqnarray}
Combining these equations we find

\begin{equation}  
$$\partial _z\partial _{\bar z}\left[ {\ln (g'^\alpha_+ f'^\alpha_- )-\ln \left(
{g_-^\alpha f_+^\alpha } \right)} \right]=-2\sum\limits_\beta  {K_{\beta \alpha }}\left(
{g'^\alpha_+ f'^\alpha_- -g_-^\alpha f_+^\alpha } \right) \ .$$  \label{74} 
\end{equation}   With the further result that

\begin{equation}  
$$\partial _{\bar z}\ln \left( {g'^\alpha_+ f_+^\alpha } \right)=\partial _z\ln \left(
{f'^\alpha_- g_-^\alpha } \right)=0 \ , $$    \end{equation}    
we recognise (\ref{74}) as the affine toda equations:

\begin{equation}  
$$\partial _z\partial _{\bar z}\phi _\alpha +\sum\limits_\beta  {K_{\beta \alpha }}\left( {\eta _+^\alpha e^{\phi _\alpha }-\eta _-^\alpha e^{-\phi _\alpha }} \right) \ .$$    \end{equation}    
Here 

\begin{equation}  
$$\phi _\alpha =\ln \left( {{{f'^\alpha_- } \over {f^\alpha_+ }}} \right)$$   
\end{equation}     and $\eta _+^\alpha$ , $\eta _-^\alpha$ are arbitrary holomorphic and
anti-holomorphic functions respectively.

\section{Some conclusions}

All of the reductions to Toda systems we have seen
use the above Result 4, although this is far from clear from a first examination.
An understanding of why this follows from the fact
that for a large class of harmonic maps, the paramatrised Lax pair can be gauge
transformed into our required form.  We need the following result.

{\bf Result 5:} {\it If the commutator $\left[ A_z,A_{\bar z} \right]$ is diagonal, then
either $A_z,A_{\bar z}$ are off-diagonal or the commutator is zero}

{\bf Proof:} Using the expansion of $A$ as in (\ref{63}), we find that since the off-diagonal elements
of the commutator are given for all $\alpha, \beta$ by

\begin{eqnarray}
h_+^\beta g_-^\alpha -h_-^\beta g_+^\alpha =0 \\
  h_+^\beta f_-^\alpha -h_-^\beta f_+^\alpha =0
\end{eqnarray}

We find that the only case when all $h^{\beta}$ do not vanish is the trivial case.

Now when we gauge transform by a group element $u$ which diagonalises the hermitian
quantity $\left[ A_z,A_{\bar z} \right]$ {\it with $u^{-1}du$ triangular}, our connection
is in the appropriate form to apply the Result 4.  In the work on
self-dual Chern-Simons theory by Dunne \cite{D}, it can be seen that for the SU(N) uniton
solutions, the necessary gauge transformation is in the required form.  These uniton solutions
are therefore Toda systems.
Moreover, Guest \cite{G}, selected {\it a priori} a connection equivalent to our form
for his reduction of the chiral model to the Toda lattice.

O'Raifeartaigh et al \cite{FWBFOR} carry out essentially this same diagonalisation process
but in a different setting.   To see this, use the decomposition in (\ref{gauss}) and 
write the commutator of the
connection of Section V of this chapter as

\begin{equation}  
$$C^{-1}B^{-1}\left[ {\tilde A_+,\tilde A_-} \right]BC$$    \end{equation}    
where

\begin{equation}  
$$\tilde A=A^{-1}dA+(dB)B^{-1}+B(dC)C^{-1}B^{-1}$$    \end{equation}    

Now we see that the quantity $BC$ corresponds to the $u$ in (\ref{u})
which diagonalises $\left[ {\tilde A_+,\tilde A_-} \right]$.
In light of the conditions (\ref{orc}), the diagonal component

\begin{equation}  
$${\left( {\partial _zx\partial _{\bar z}y-\partial _zx\partial _{\bar z}y} \right)e^{-\phi }}$$    \end{equation}    

is just the Liouville field $e^\phi $ up 
to a holomorphic function.  This corresponds
to the SU(1,1) Toda solution
 result in (\ref{md}) and confirms the solution found by Fujii \cite{F}.

\end{document}